\documentclass[11pt]{article}
\usepackage[english]{babel}
\usepackage[cp1251]{inputenc}
\usepackage[psamsfonts]{amssymb}
\usepackage[intlimits]{amsmath}
\usepackage{amsfonts}

\begin{document}

\begin{center}
{\bf ERROR ESTIMATES FOR THE ITERATIVELY REGULARIZED GRADIENT METHOD FOR ILL-POSED OPTIMIZATION PROBLEMS IN A HILBERT SPACE}\footnote{This work was supported by the Russian Science Foundation (project no. 22-71-10070)}
\\[2ex]
{\bf M.\,M.\,Kokurin}\\[3ex]
\end{center}

\medskip

We study the iteratively regularized gradient method applied to the ill-posed problem of minimizing a general smooth functional on a Hilbert space. Under a source condition on the solution, we establish error estimates for this method when the objective functional is known exactly or only approximately. For these estimates we provide two proofs based on substantially different ideas and compare the assumptions and conclusions of the corresponding theorems. In addition, we present an improved proof of an error estimate for the Tikhonov method applied to ill-posed optimization problems with an exactly known functional. We also prove a theorem on the existence of a global minimizer of the Tikhonov functional without assuming weak lower semicontinuity of the original functional.

\medskip


\section{Problem statement}

Consider the minimization problem
$$
J(x)\to\min\limits_{x\in H} \eqno{(1)}
$$
for a nonlinear functional $J:\,H\to\mathbb{R}$ on a real Hilbert space $H$. We assume throughout that problem (1) has a global minimizer $x^*\in H$. We assume that $J$ is twice continuously Fr\'echet differentiable and that its second Fr\'echet derivative satisfies the Lipschitz condition
$$
\|J^{\prime\prime}(x)-J^{\prime\prime}(y)\|_{L(H)}\leq L\|x-y\|,\quad x,y\in H \eqno{(2)}
$$
for some constant $L>0$. Here and below, $\|\cdot\|$ denotes the norm in $H$, while $\|\cdot\|_{L(H)}$ denotes the operator norm on $L(H)$, the space of bounded linear operators on $H$.

In general, optimization problem (1) is ill-posed (see, e.g., [1, 2]). This means that it cannot be solved directly by classical minimization methods: even if such a method produces a minimizing sequence $\{x_n\}\subset H$ such that $\lim\limits_{n\to\infty}J(x_n)=J(x^*)$, this sequence need not converge to the sought solution $x^*$ in the norm of $H$. Regularization methods are used to solve ill-posed optimization problems. The best-known of these methods is the Tikhonov method. One of its simplest versions replaces problem (1) by the problem of minimizing the Tikhonov functional
$$
T_\alpha:\,H\to\mathbb{R},\quad T_\alpha(x)=J(x)+\alpha\|x-\xi\|^2. \eqno{(3)}
$$
Here $\alpha>0$ is a regularization parameter, and $\xi\in H$ is specified a priori and may be regarded as an approximation to the sought minimizer $x^*$. It is known that if the functional $J$ is weakly lower semicontinuous, then the functional $T_\alpha$ has a global minimizer $u_\alpha\in H$. In the Tikhonov method, such a minimizer is used as an approximation to $x^*$. If, in addition, $x^*$ is the unique global minimizer of $J$ closest to $\xi$, then $u_\alpha\to x^*$ in the norm of $H$ as $\alpha\to 0$. In general, however, $u_\alpha$ need not exist; even when it does, finding it exactly may be difficult. In that case, one instead seeks a point $u_{\alpha,\varepsilon}\in H$, $\varepsilon>0$, satisfying the two-sided inequality
$$
\inf\limits_{x\in H}T_\alpha(x)\leq T_\alpha(u_{\alpha,\varepsilon})\leq \inf\limits_{x\in H}T_\alpha(x)+\varepsilon, \eqno{(4)}
$$
where $\varepsilon$ must be chosen in coordination with $\alpha$ (see [1, Chap.~2, \S~5]).

Unlike the original functional $J$, the Tikhonov functional (3) can often be minimized directly by classical iterative methods. The simplest choice is the gradient method:
$$
x_{n+1}=x_n-\mu T_\alpha^\prime(x_n),\quad n\in\mathbb{N},\quad\mu>0. \eqno{(5)}
$$
Here
$$
T_\alpha^\prime(x)=J^\prime(x)+2\alpha(x-\xi). \eqno{(6)}
$$
Method (5) can be applied directly as follows. First, one chooses a sufficiently small value of $\alpha$ and performs several gradient iterations to approach $u_\alpha$ or to obtain a point satisfying (4) for a prescribed value of $\varepsilon$. One then decreases $\alpha$ and repeats the procedure, possibly with a different value of $\mu$. This raises the questions of how many gradient iterations to perform and how to choose $\mu$ for each $\alpha$. The one-step iteratively regularized gradient method is particularly convenient: exactly one step of (5) is performed for each value of $\alpha$. More precisely, we choose positive sequences $\{\alpha_n\}$ and $\{\mu_n\}$, with $\{\alpha_n\}\to 0$, together with an initial approximation $x_1$, and define
$$
x_{n+1}=x_n-\mu_n T_{\alpha_n}^\prime(x_n),\quad n\in\mathbb{N}. \eqno{(7)}
$$
The resulting scheme is the iteratively regularized gradient method.

An additional difficulty arises when, as is often the case in applications, the functional $J$ is not known exactly. In this paper, we assume that an approximation $J_\delta:\,H\to\mathbb{R}$ is available in place of the exact functional $J$ and satisfies
$$
\|J_\delta^\prime(x)-J^\prime(x)\|\leq\delta(1+\|x-\xi\|),\quad x\in H \eqno{(8)}
$$
with a known error level $\delta\geq 0$. Condition (8) on the approximate functional was used, for example, in [3, 4]. Some of the results below require a stronger condition on $J_\delta$ in addition to (8). When the functional $J$ is known only approximately, the Tikhonov method uses the approximate Tikhonov functional
$$
T_\alpha^\delta:\,H\to\mathbb{R},\quad T_\alpha^\delta(x)=J_\delta(x)+\alpha\|x-\xi\|^2,
$$
in place of (3), and the iteratively regularized gradient method uses the iteration
$$
x_{n+1}=x_n-\mu_n T_{\alpha_n}^{\delta\,\prime}(x_n)=x_n-\mu_n\bigl(J_\delta^\prime(x_n)+2\alpha_n(x_n-\xi)\bigr),\quad n\in\mathbb{N}. \eqno{(9)}
$$
Below, we regard the case of an exactly known functional $J$ as a special case of the general problem setting with $J_\delta=J$, $\delta=0$. In this case, scheme (9) reduces to (7).

For a nonzero error level $\delta>0$, iterative regularization methods such as (9) need not converge to the desired solution as $n\to\infty$. They must therefore be stopped at an index $N=N(\delta)$, and $x_{N(\delta)}$ is taken as the final approximation to $x^*$. In the exact case $\delta=0$, no stopping rule is needed.

This leads to the following questions. Under what conditions and at what rate does the sequence generated by (7) converge to $x^*$? Given a suitable stopping rule $N=N(\delta)$, under what conditions and at what rate does the stopped iterate $x_{N(\delta)}$ generated by (9) converge to $x^*$ as $\delta\to 0$? In both cases, convergence is understood in the norm of $H$. Since iterations (7) and (9) are derived from the gradient method applied to the Tikhonov functional, these questions are also closely related to error estimates for the Tikhonov method. Note that error estimates for methods of solving ill-posed optimization problems can usually be obtained only under restrictive a priori assumptions on the sought solution. In this paper, we use the source condition
$$
\exists w\in H:\,x^*-\xi=J^{\prime\prime}(x^*)w \eqno{(10)}
$$
together with additional bounds on $\|w\|$.

We now give a brief overview of known results related to these questions. The iteratively regularized gradient method for optimization problems was introduced in [5], but only for the minimization of convex functionals. In [6], the iteratively regularized gradient method was studied for nonlinear operator equations $F(x)=0$, where $F$ is a general smooth operator between Hilbert spaces. Error estimates $O(\sqrt{\alpha_n})$ and $O(\sqrt\delta)$ were obtained in the exact and perturbed cases, respectively. For an approximately specified $F$, an a priori stopping rule was used. This line of research was further developed in [7], where an error estimate $O(\sqrt\delta)$ was obtained for the first time for the iteratively regularized gradient method with an a posteriori stopping rule. However, this result again concerned nonlinear operator equations with an error of level $\delta$ in the right-hand side rather than general optimization problems. Nevertheless, [7] imposed none of the restrictive structural assumptions on the operator $F$ that are used in many related works (e.g., [3]); only sufficient smoothness was required.

We also mention [8], which studied the convergence of the gradient method without a Tikhonov term for the problem of minimizing a nonconvex functional in a Hilbert space under a structural assumption on the nonlinearity of the functional. This line of research was continued in [9], where conditions for the strong and weak convergence of several methods for constrained minimization of functionals in Hilbert and Banach spaces, including the projected gradient method, were studied. Applications of these methods to specific inverse problems were also discussed. For a detailed survey of the general theory of ill-posed problems, including ill-posed optimization problems and methods for solving them, see [2].

Convergence theorems for the Tikhonov method and several other basic methods for solving ill-posed optimization problems, including constrained problems, are presented in [1]. In [10], power-type error estimates were first obtained for the Tikhonov method applied to a general smooth nonconvex functional $J$ on a Hilbert space: $O(\sqrt\alpha)$ when $J$ is known exactly and $O(\sqrt[4]{\delta})$ when it is known approximately. In [4], these estimates were improved to $O(\alpha)$ and $O(\sqrt\delta)$, respectively.

The present paper is a natural continuation of [4]. Here, for the first time, we prove power-type error estimates for the iteratively regularized gradient method applied to the problem of minimizing a general smooth functional on a Hilbert space, without imposing structural assumptions on the nonlinearity. These estimates are $O(\alpha_n)$ for method (7) in the exact case and $O(\sqrt\delta)$ for the stopped iteration (9) in the perturbed case. Here $\alpha_n=\alpha_0(n+k)^{-s}$ and $\mu_n=\mu\alpha_n$, $n\in\mathbb{N}$, where $\alpha_0$, $k$, $s$, and $\mu>0$ are suitable parameters. We give two substantially different proofs of these estimates. The first proof follows the strategy outlined in [6] and [11, Chap.~2, \S~8]: we first estimate the distance between the iterates and the corresponding minimizers of the Tikhonov functional and then combine this bound with an error estimate for the Tikhonov method itself. A similar approach is described in [12, Chap.~6], for the construction and analysis of two-stage iterative methods for solving operator equations. As part of this approach, in Section 2 we give, for the exact case $\delta=0$, a proof of a Tikhonov error estimate that is substantially shorter than the proof given in [4]. In Section 3, we prove two auxiliary results that are also of independent interest: the existence of a global minimizer of the Tikhonov functional without assuming weak lower semicontinuity of $J$, and an estimate of the distance between such minimizers for different values of the regularization parameter. The first proof of the main result is presented in Section 4. The second proof, given in Section 5, estimates the error of the iteratively regularized gradient method directly, without using any properties of the Tikhonov functional or its minimizers. Both proofs apply to problems with arbitrary functionals $J$ satisfying condition (2), provided that the element $w$ in source condition (10) has sufficiently small norm; precise conditions on this norm are stated explicitly. Finally, we compare the two approaches and consider the case in which the functional is only locally smooth.

\section{An error estimate for the Tikhonov method in the case $\delta=0$ with an explicit formula for the coefficient}

This section provides an auxiliary result. In the first proof of the error estimate for the iteratively regularized gradient method (9), we will rely on the following error estimate for the Tikhonov method:
$$
\|u_{\alpha,\varepsilon}-x^*\|\leq C_0\alpha, \eqno{(11)}
$$
which holds when the parameters $\varepsilon$ and $\alpha$ are chosen in a coordinated manner. Here and below, $C_0,\,C_1,\,C_2,\,\ldots$ denote nonnegative constants that are independent of $\alpha$, $\varepsilon$, and $\delta$, but may depend, for example, on the exact functional $J$ and the choice of the point $\xi$. Estimate (11) was obtained in [4] under conditions (2), (10), and the additional assumption
$$
L\|w\|<1, \eqno{(12)}
$$
with the choice $\varepsilon=K\alpha^4$, $K>0$. However, to derive an error estimate for the iteratively regularized gradient method, it is important to know not only that a constant $C_0$ satisfying (11) exists, but also how this constant depends on $L$ and $\|w\|$. This dependence was not stated explicitly in [4]. Therefore, in this section we present a proof of (11) that yields a simple explicit formula for $C_0$. This proof is also substantially shorter than the one given in [4] and is valid under a slightly weaker condition on $\varepsilon$ (see formula (13) below). For our purposes, an explicit formula for the coefficient $C_0$ is needed only in the case $\delta=0$. By contrast, [4] treated the general case $\delta\geq 0$ from the outset.

Accordingly, let $u=u_{\alpha,\varepsilon}$ satisfy (4). If $\varepsilon=0$, then $u_{\alpha,0}=u_\alpha$, provided that this element exists, whereas for $\varepsilon>0$ the existence of elements $u=u_{\alpha,\varepsilon}$ is guaranteed.

\medskip

{\bf Theorem 2.1.} Let the functional $J:\,H\to\mathbb{R}$ on a real Hilbert space $H$ have a global minimizer $x^*$ and satisfy the smoothness condition (2). Suppose, in addition, that source condition (10), inequality (12), and the relation
$$
0\leq\varepsilon\leq K\alpha^3 \eqno{(13)}
$$
hold with some fixed constant $K\geq 0$. Then estimate (11) holds for the Tikhonov method with the exact functional $J$.

{\bf Proof.} For convenience, the proof below is divided into two stages. The first stage largely follows the argument in [4], with some simplifications. Its key idea is to derive estimate (14) and then substitute $z=x^*-2\alpha w$. This substitution can be motivated by the following heuristic argument. If $\varepsilon=0$, then $T_\alpha^\prime(u_\alpha)=0_H$, so $J^\prime(u_\alpha)=-2\alpha(u_\alpha-\xi)=-2\alpha(u_\alpha-x^*)-2\alpha(x^*-\xi)$. On the other hand, $J^\prime(u_\alpha)=J^\prime(u_\alpha)-J^\prime(x^*)\approx J^{\prime\prime}(x^*)(u_\alpha-x^*)$. Therefore, $(J^{\prime\prime}(x^*)+2\alpha E)(u_\alpha-x^*)\approx -2\alpha(x^*-\xi)=-2\alpha J^{\prime\prime}(x^*)w$. Hence, $u_\alpha-x^*\approx -2\alpha(J^{\prime\prime}(x^*)+2\alpha E)^{-1}J^{\prime\prime}(x^*)w\approx -2\alpha w$. Here and below, $0_H$ is the zero element of $H$, and $E:\,H\to H$ is the identity operator. The last approximate equality is based on the fact that the element $w$ in condition (10) can always be chosen in $({\rm{Ker}}\,J^{\prime\prime}(x^*))^\perp$, in which case $(J^{\prime\prime}(x^*)+2\alpha E)^{-1}J^{\prime\prime}(x^*)w\,\to\,w$ in the norm of $H$ as $\alpha\to 0$. Thus, $x^*-2\alpha w$ is a natural approximation to $u_\alpha$ and, for sufficiently small $\alpha$ and $\varepsilon$, also to $u=u_{\alpha,\varepsilon}$.

In contrast, the second stage of the proof differs substantially from the argument in [4]. In [4], the difference $J(x^*-2\alpha w)-J(u)$ was estimated using an auxiliary self-adjoint operator $B$ such that $B(x^*-\xi)=w$. Instead, below we use the inequality $J(x^*)\leq J(u+2\alpha w)$, where $u+2\alpha w$ is chosen to match the substitution $z=x^*-2\alpha w$. This makes the proof substantially simpler and shorter.

\medskip

{\bf Stage 1.} It follows from (4) that, for any $z\in H$,
$$
J(u)+\alpha\|u-\xi\|^2\leq J(z)+\alpha\|z-\xi\|^2+\varepsilon,
$$
and hence
$$
\|u-\xi\|^2-\|z-\xi\|^2\leq\frac{J(z)-J(u)+\varepsilon}{\alpha}.
$$
Using the readily verified identity
$$
\|u-\xi\|^2-\|z-\xi\|^2=\|u-x^*\|^2-2(x^*-\xi,z-u)-\|z-x^*\|^2
$$
and source condition (10), we obtain
$$
\|u-x^*\|^2\leq\frac{J(z)-J(u)}{\alpha}+2(J^{\prime\prime}(x^*)w,z-u)+\|z-x^*\|^2+\frac{\varepsilon}{\alpha}. \eqno{(14)}
$$
Substituting $z=x^*-2\alpha w$ into this inequality gives
\begin{equation*}
\begin{split}
&\|u-x^*\|^2\leq\frac{J(x^*-2\alpha w)-J(u)}{\alpha}+2(J^{\prime\prime}(x^*)w,x^*-u)-\\
&-4\alpha(J^{\prime\prime}(x^*)w,w)+4\alpha^2\|w\|^2+\frac{\varepsilon}{\alpha}.
\end{split}
\tag{15}
\end{equation*}

Our immediate objective is to obtain a convenient representation of the term $(J^{\prime\prime}(x^*)w,x^*-u)$ on the right-hand side of (15). To this end, observe that
$$
J^\prime(u)=J^\prime(u)-J^\prime(x^*)=\int\limits_0^1 J^{\prime\prime}(x^*+t(u-x^*))(u-x^*)dt=J^{\prime\prime}(x^*)(u-x^*)+I,
$$
where
$$
I=\int\limits_0^1\bigl(J^{\prime\prime}(x^*+t(u-x^*))-J^{\prime\prime}(x^*)\bigr)(u-x^*)dt,\quad \|I\|\leq\frac{L}{2}\|u-x^*\|^2. \eqno{(16)}
$$
It follows that
$$
(T_\alpha^\prime(u),w)=(J^\prime(u),w)+2\alpha(u-\xi,w)=
$$
$$
=(J^{\prime\prime}(x^*)(u-x^*),w)+(I,w)+2\alpha(u-\xi,w).
$$
Using the self-adjointness of $J^{\prime\prime}(x^*)$, we obtain
$$
(J^{\prime\prime}(x^*)w,x^*-u)=-(J^{\prime\prime}(x^*)(u-x^*),w)=(I,w)+2\alpha(u-\xi,w)-(T_\alpha^\prime(u),w). \eqno{(17)}
$$

We now substitute representation (17) into (15):
\begin{equation*}
\begin{split}
&\|u-x^*\|^2\leq\frac{J(x^*-2\alpha w)-J(u)}{\alpha}+2(I,w)+4\alpha(u-\xi,w)-\\
&-2(T_\alpha^\prime(u),w)-4\alpha(J^{\prime\prime}(x^*)w,w)+4\alpha^2\|w\|^2+\frac{\varepsilon}{\alpha}.
\end{split}
\tag{18}
\end{equation*}
Here
$$
4\alpha(u-\xi,w)=4\alpha(u-x^*,w)+4\alpha(x^*-\xi,w)=4\alpha(u-x^*,w)+4\alpha(J^{\prime\prime}(x^*)w,w).
$$
Upon substituting this identity into (18), the two terms $4\alpha(J^{\prime\prime}(x^*)w,w)$ with opposite signs on its right-hand side cancel:
$$
\|u-x^*\|^2\leq \frac{J(x^*-2\alpha w)-J(u)}{\alpha}+2(I,w)+
$$
$$
+4\alpha(u-x^*,w)-2(T_\alpha^\prime(u),w)+4\alpha^2\|w\|^2+\frac{\varepsilon}{\alpha}.
$$
Using (16), we estimate the term $2(I,w)$ on the right-hand side by $L\|w\|\|u-x^*\|^2$ and move it to the left-hand side:
$$
(1-L\|w\|)\|u-x^*\|^2\leq \frac{J(x^*-2\alpha w)-J(u)}{\alpha}+
$$
$$
+4\alpha(u-x^*,w)-2(T_\alpha^\prime(u),w)+4\alpha^2\|w\|^2+\frac{\varepsilon}{\alpha}.
$$
We now divide both sides of the resulting inequality by $1-L\|w\|$. This is legitimate by (12). We arrive at the intermediate estimate
\begin{equation*}
\begin{split}
\|u-x^*\|^2\leq\frac{1}{1-L\|w\|}\Bigl( &\frac{J(x^*-2\alpha w)-J(u)}{\alpha}+4\alpha(u-x^*,w)-\\
&-2(T_\alpha^\prime(u),w)+4\alpha^2\|w\|^2+\frac{\varepsilon}{\alpha} \Bigr).
\end{split}
\tag{19}
\end{equation*}
This estimate will serve as the basis for the next stage of the proof.

\medskip

{\bf Stage 2.} We now derive an estimate for the expression $J(x^*-2\alpha w)-J(u)$ in (19). We add and subtract $J(x^*)$ in this expression and then use the inequality $J(x^*)\leq J(u+2\alpha w)$:
\begin{equation*}
\begin{split}
&J(x^*-2\alpha w)-J(u)=J(x^*-2\alpha w)-J(x^*)+J(x^*)-J(u)\leq\\
&\leq \bigl(J(x^*-2\alpha w)-J(x^*)\bigr)+\bigl(J(u+2\alpha w)-J(u)\bigr).
\end{split}
\tag{20}
\end{equation*}
By the Taylor estimate implied by (2),
$$
J(x+h)\leq J(x)+(J^\prime(x),h)+\frac{1}{2}(J^{\prime\prime}(x)h,h)+\frac{1}{6}L\|h\|^3,\quad x,h\in H,
$$
the first parenthesized expression on the right-hand side of (20) satisfies
\begin{equation*}
\begin{split}
&J(x^*-2\alpha w)-J(x^*)\leq (J^\prime(x^*),-2\alpha w)+2\alpha^2(J^{\prime\prime}(x^*)w,w)+\frac{4}{3}L\|w\|^3\alpha^3=\\
&=2\alpha^2(J^{\prime\prime}(x^*)w,w)+\frac{4}{3}L\|w\|^3\alpha^3.
\end{split}
\tag{21}
\end{equation*}
Here we have used the fact that $J^\prime(x^*)=0_H$, since $x^*$ is a minimizer of $J$.

Similarly, the second parenthesized expression on the right-hand side of (20) satisfies
$$
J(u+2\alpha w)-J(u)\leq 2\alpha (J^\prime(u),w)+2\alpha^2(J^{\prime\prime}(u)w,w)+\frac{4}{3}L\|w\|^3\alpha^3.
$$
By (6),
$$
(J^\prime(u),w)=(T_\alpha^\prime(u),w)-2\alpha(u-\xi,w)=
$$
$$
=(T_\alpha^\prime(u),w)-2\alpha(u-x^*,w)-2\alpha(J^{\prime\prime}(x^*)w,w).
$$
Furthermore,
$$
(J^{\prime\prime}(u)w,w)=(J^{\prime\prime}(x^*)w,w)+\bigl((J^{\prime\prime}(u)-J^{\prime\prime}(x^*))w,w\bigr)\leq
$$
$$
\leq(J^{\prime\prime}(x^*)w,w)+L\|u-x^*\|\|w\|^2,
$$
and therefore
\begin{equation*}
\begin{split}
&J(u+2\alpha w)-J(u)\leq 2\alpha(T_\alpha^\prime(u),w)-4\alpha^2(u-x^*,w)-4\alpha^2(J^{\prime\prime}(x^*)w,w)+\\
&+2\alpha^2(J^{\prime\prime}(x^*)w,w)+2L\|w\|^2\alpha^2\|u-x^*\|+\frac{4}{3}L\|w\|^3\alpha^3.
\end{split}
\tag{22}
\end{equation*}

Substituting (21) and (22) into (20), we obtain
$$
J(x^*-2\alpha w)-J(u)\leq (2\alpha^2-4\alpha^2+2\alpha^2)(J^{\prime\prime}(x^*)w,w)+\frac{4}{3}L\|w\|^3\alpha^3+
$$
$$
+2\alpha(T_\alpha^\prime(u),w)-4\alpha^2(u-x^*,w)+2L\|w\|^2\alpha^2\|u-x^*\|+\frac{4}{3}L\|w\|^3\alpha^3=
$$
$$
=2\alpha(T_\alpha^\prime(u),w)-4\alpha^2(u-x^*,w)+2L\|w\|^2\alpha^2\|u-x^*\|+\frac{8}{3}L\|w\|^3\alpha^3.
$$
Note that the cancellation of like terms proportional to $(J^{\prime\prime}(x^*)w,w)$ here is a consequence of the coordinated choice of the points $x^*-2\alpha w$ and $u+2\alpha w$ in (20). Substitution of the resulting inequality into (19) also cancels the terms $2(T_\alpha^\prime(u),w)$ and $4\alpha(u-x^*,w)$:
$$
\|u-x^*\|^2\leq\frac{1}{1-L\|w\|}\left( 2L\|w\|^2\alpha\|u-x^*\|+\left(4\|w\|^2+\frac{8}{3}L\|w\|^3\right)\alpha^2+\frac{\varepsilon}{\alpha} \right).
$$
Using (13), we arrive at
$$
\|u-x^*\|^2\leq C_1\alpha\|u-x^*\|+C_2(K)\alpha^2,
$$
where
$$
C_1=\frac{2L\|w\|^2}{1-L\|w\|},\quad C_2(K)=\frac{4\|w\|^2+\frac{8}{3}L\|w\|^3+K}{1-L\|w\|}.
$$
In other words, the quantity $y=\|u-x^*\|/\alpha\geq 0$ satisfies the quadratic inequality
$$
y^2-C_1y-C_2(K)\leq 0
$$
and therefore does not exceed its positive root $(C_1+\sqrt{C_1^2+4C_2(K)})/2$. Thus,
$$
\|u-x^*\|\leq C_0(K)\alpha,\quad C_0(K)=\frac{C_1+\sqrt{C_1^2+4C_2(K)}}{2}. \eqno{(23)}
$$
{\bf Theorem proved.}

\medskip

We have thus obtained an error estimate, with an explicit coefficient $C_0(K)$, for the approximation $u=u_{\alpha,\varepsilon}$ produced by the Tikhonov method applied to problem (1) with the exact functional $J$. This estimate will be used below to derive an error estimate for the iteratively regularized gradient method. To do so, however, we need additional results on the existence and properties of exact minimizers of the Tikhonov functional, which are also of independent interest. The next section is devoted to these results.

\section{Exact minimizers of the Tikhonov functional and their properties}

Define the constant
$$
\widehat C_0=C_0(0)=\frac{1}{1-L\|w\|}\Biggl(L\|w\|^2+
$$
$$
+\sqrt{L^2\|w\|^4+(1-L\|w\|)\left(4\|w\|^2+\frac{8}{3}L\|w\|^3\right)}\Biggr),
$$
which is the coefficient in formula (23) for $K=0$. In addition to (12), we will assume below that
$$
L\widehat C_0<2. \eqno{(24)}
$$
In Section 4, condition (24) will be strengthened to (42), and an exact criterion for the validity of (42) will be derived.

We prove that, under condition (24), the Tikhonov functional (3) has a global minimizer on the entire space $H$, without assuming weak lower semicontinuity of $J$.

\medskip

{\bf Theorem 3.1.} Suppose that conditions (2), (10), (12), and (24) hold. Then, for every $\alpha>0$, the functional $T_\alpha$ has a unique global minimizer $u_\alpha$, and
$$
\|u_\alpha-x^*\|\leq \widehat C_0\alpha. \eqno{(25)}
$$

{\bf Proof.} By the continuity of $C_0(K)$ and the strict inequality (24), there exists a sufficiently small $K>0$ such that $LC_0(K)<2$. Fix $\alpha>0$ and choose a sequence $\{\varepsilon_j\}\to 0$ contained in the interval $(0,K\alpha^3]$. By Theorem 2.1, all the points $v_j=u_{\alpha,\varepsilon_j}$ belong to the ball
$$
B_\alpha=\{x\in H\,|\, \|x-x^*\|\leq C_0(K)\alpha\}.
$$
Since $x^*$ is a minimizer of the twice continuously differentiable functional $J$, the operator $J^{\prime\prime}(x^*)$ is positive semidefinite. Therefore, for arbitrary $x\in B_\alpha$ and $h\in H$,
$$
(T_\alpha^{\prime\prime}(x)h,h)=((J^{\prime\prime}(x)+2\alpha E)h,h)=
$$
$$
=(J^{\prime\prime}(x^*)h,h)+\bigl((J^{\prime\prime}(x)-J^{\prime\prime}(x^*))h,h\bigr)+2\alpha\|h\|^2\geq
$$
$$
\geq (2\alpha-\|J^{\prime\prime}(x)-J^{\prime\prime}(x^*)\|_{L(H)})\|h\|^2\geq (2-LC_0(K))\alpha\|h\|^2.
$$
Since $LC_0(K)<2$, the functional $T_\alpha$ is strongly convex on the ball $B_\alpha$. Set $m_\alpha=\inf\limits_{x\in H}T_\alpha(x)$. The strong convexity of $T_\alpha$ on the line segment joining $v_j$ and $v_l$ gives
$$
m_\alpha\leq T_\alpha\left(\frac{v_j+v_l}{2}\right)\leq\frac{T_\alpha(v_j)+T_\alpha(v_l)}{2}-\frac{(2-LC_0(K))\alpha}{8}\|v_j-v_l\|^2\leq
$$
$$
\leq\frac{2m_\alpha+\varepsilon_j+\varepsilon_l}{2}-\frac{(2-LC_0(K))\alpha}{8}\|v_j-v_l\|^2.
$$
Consequently,
$$
\|v_j-v_l\|^2\leq\frac{4(\varepsilon_j+\varepsilon_l)}{(2-LC_0(K))\alpha}.
$$
Thus, $\{v_j\}$ is a Cauchy sequence and converges in $H$. By construction, it is a minimizing sequence for $T_\alpha$ on the entire space $H$, and the continuity of $T_\alpha$ implies that its limit is a global minimizer of this functional. This minimizer is unique because any other global minimizer belongs to the ball $B_\alpha$ by Theorem 2.1 with $\varepsilon=0$ and the chosen $K$, and $T_\alpha$ is strongly convex on this ball. We emphasize that $u_\alpha$ minimizes $T_\alpha$ on the entire space $H$, not merely on the ball $B_\alpha$. Estimate (25) follows directly from Theorem 2.1 with $K=0$. {\bf Theorem proved.}

\medskip

The following lemma gives an estimate for the distance between minimizers of the Tikhonov functional corresponding to different values of the regularization parameter. An analogous estimate was used in [6] to prove an error estimate for the iteratively regularized gradient method applied to operator equations. In Section 4, we adapt this argument to the general optimization problem (1).

\medskip

{\bf Lemma 3.1.} Let $0<\beta<\alpha\leq\bar\alpha$, where $\bar\alpha>0$ is fixed, and suppose that the assumptions of Theorem 3.1 hold. Then
$$
\|u_\alpha-u_\beta\|\leq C_3\frac{\alpha-\beta}{\sqrt\alpha}. \eqno{(26)}
$$

{\bf Proof.} We have $T_\alpha^\prime(u_\alpha)=T_\beta^\prime(u_\beta)=0_H$, i.e.,
$$
J^\prime(u_\alpha)+2\alpha(u_\alpha-\xi)=0_H,\quad J^\prime(u_\beta)+2\beta(u_\beta-\xi)=0_H.
$$
Hence,
$$
J^\prime(u_\alpha)-J^\prime(u_\beta)=-2\alpha(u_\alpha-\xi)+2\beta(u_\beta-\xi)=-2\alpha(u_\alpha-u_\beta)-2(\alpha-\beta)(u_\beta-\xi).
$$
Taking the inner product of both sides with $u_\alpha-u_\beta$ and expanding $J^\prime(u_\alpha)$ about $u_\beta$ by Taylor's formula with an integral remainder, we obtain
$$
2\alpha\|u_\alpha-u_\beta\|^2=-\int\limits_0^1\bigl(J^{\prime\prime}(u_\beta+t(u_\alpha-u_\beta))(u_\alpha-u_\beta),
u_\alpha-u_\beta\bigr)dt-
$$
$$
-2(\alpha-\beta)(u_\beta-\xi,u_\alpha-u_\beta).
$$
We isolate the term $J^{\prime\prime}(x^*)$ in the integral and use the identity $u_\beta-\xi=(u_\beta-x^*)+J^{\prime\prime}(x^*)w$, which follows from source condition (10):
$$
2\alpha\|u_\alpha-u_\beta\|^2=\int\limits_0^1\Bigl(\bigl(J^{\prime\prime}(x^*)-J^{\prime\prime}(u_\beta+t(u_\alpha-u_\beta))\bigr)(u_\alpha-u_\beta),
u_\alpha-u_\beta\Bigr)dt-
$$
$$
-\bigl(J^{\prime\prime}(x^*)(u_\alpha-u_\beta),u_\alpha-u_\beta\bigr)-2(\alpha-\beta)((u_\beta-x^*)+J^{\prime\prime}(x^*)w,u_\alpha-u_\beta).
$$
The entire line segment joining $u_\beta$ and $u_\alpha$ lies in the ball $\|x-x^*\|\leq \widehat C_0\alpha$, and $\|u_\beta-x^*\|\leq \widehat C_0\beta$. Estimating the norm of the difference of the second derivatives of $J$ using (2), we arrive at
\begin{equation*}
\begin{split}
&(2-L\widehat C_0)\alpha\|u_\alpha-u_\beta\|^2\leq
-(J^{\prime\prime}(x^*)(u_\alpha-u_\beta),u_\alpha-u_\beta)-\\
&-2(\alpha-\beta)(J^{\prime\prime}(x^*)w,u_\alpha-u_\beta)
+2\widehat C_0\beta(\alpha-\beta)\|u_\alpha-u_\beta\|.
\end{split}
\end{equation*}
The first two terms on the right-hand side are estimated by completing the square:
\begin{equation*}
\begin{split}
&-(J^{\prime\prime}(x^*)(u_\alpha-u_\beta),u_\alpha-u_\beta)-2(\alpha-\beta)(J^{\prime\prime}(x^*)w,u_\alpha-u_\beta)=\\
&=(\alpha-\beta)^2\left\|\bigl(J^{\prime\prime}(x^*)\bigr)^{1/2}w\right\|^2-\left\|\bigl(J^{\prime\prime}(x^*)\bigr)^{1/2}\bigl(u_\alpha-u_\beta+(\alpha-\beta)w\bigr)\right\|^2\leq\\
&\leq(\alpha-\beta)^2\left\|\bigl(J^{\prime\prime}(x^*)\bigr)^{1/2}w\right\|^2.
\end{split}
\end{equation*}
To estimate the remaining term, we use the inequality
$$
2\widehat C_0\beta(\alpha-\beta)\|u_\alpha-u_\beta\|\leq \eta\|u_\alpha-u_\beta\|^2+\frac{\widehat C_0^2\beta^2(\alpha-\beta)^2}{\eta},\quad\eta=\frac{(2-L\widehat C_0)\alpha}{2}.
$$
We obtain
$$
\|u_\alpha-u_\beta\|^2\leq\frac{2\left\|\bigl(J^{\prime\prime}(x^*)\bigr)^{1/2}w\right\|^2}{2-L\widehat C_0}\frac{(\alpha-\beta)^2}{\alpha}+\frac{4\widehat C_0^2}{(2-L\widehat C_0)^2}
\frac{\beta^2}{\alpha}\cdot\frac{(\alpha-\beta)^2}{\alpha}.
$$
Since $\beta^2/\alpha\leq \alpha\leq\bar\alpha$, estimate (26) follows directly with the constant
$$
C_3=\left(\frac{2\left\|\bigl(J^{\prime\prime}(x^*)\bigr)^{1/2}w\right\|^2}{2-L\widehat C_0}+\frac{4\widehat C_0^2\bar\alpha}{(2-L\widehat C_0)^2}\right)^{\frac{1}{2}}.
$$
{\bf Lemma proved.}

\medskip

We now derive an error estimate for the iteratively regularized gradient method.

\section{An error estimate for the iteratively regularized gradient method: first proof}

In this section, in addition to (8), we assume that the approximate functional $J_\delta$ is twice continuously Fr\'echet differentiable and satisfies
$$
\|J_\delta^{\prime\prime}(x)-J^{\prime\prime}(x)\|_{L(H)}\leq\delta(1+\|x-\xi\|),\quad x\in H. \eqno{(27)}
$$
If $\delta=0$, both conditions (8) and (27) are satisfied by the exact functional $J_\delta=J$. Thus, the arguments below cover both the exact and the perturbed cases: in the former, $J_\delta=J$ and $\delta=0$, whereas in the latter, $\delta>0$.

We first consider one step of the iteratively regularized gradient method (9) with $\alpha_n=\alpha$ and $\mu_n=\mu\alpha$, where $\alpha,\mu>0$:
$$
x^+=x-\mu\alpha T_\alpha^{\delta\,\prime}(x)=x-\mu\alpha\bigl(J_\delta^\prime(x)+2\alpha(x-\xi)\bigr). \eqno{(28)}
$$
As in Lemma 3.1, we fix $\bar\alpha>0$ and assume that
$$
0<\alpha\leq\bar\alpha. \eqno{(29)}
$$
Since step (28) is obtained by applying one iteration of the standard gradient method to the approximate Tikhonov functional $T_\alpha^\delta$, it is natural to expect that, for a sufficiently small error level $\delta>0$ and under suitable conditions, the point $x^+$ will be closer to the minimizer $u_\alpha$ of the exact functional $T_\alpha$ than the original point $x$. Our aim is to estimate $\|x^+-u_\alpha\|$ in terms of the initial distance $\|x-u_\alpha\|$. We assume that
$$
\|x-u_\alpha\|\leq \varkappa\alpha, \eqno{(30)}
$$
where restrictions on the coefficient $\varkappa>0$ will be imposed below. We also assume that the hypotheses of Theorem 3.1 hold.

Using the identity
$$
T_\alpha^{\delta\,\prime}(u_\alpha)=J_\delta^\prime(u_\alpha)+2\alpha(u_\alpha-\xi)=\bigl(J_\delta^\prime(u_\alpha)-J^\prime(u_\alpha)\bigr)+T_\alpha^\prime(u_\alpha)=J_\delta^\prime(u_\alpha)-J^\prime(u_\alpha),
$$
we write
$$
x^+-u_\alpha=x-u_\alpha-\mu\alpha \bigl( (T_\alpha^{\delta\,\prime}(x)-T_\alpha^{\delta\,\prime}(u_\alpha))+T_\alpha^{\delta\,\prime}(u_\alpha) \bigr)=
$$
$$
=x-u_\alpha-\mu\alpha \left( \int\limits_0^1 T_\alpha^{\delta\,\prime\prime}
(u_\alpha+t(x-u_\alpha))(x-u_\alpha)dt+\bigl(J_\delta^\prime(u_\alpha)-J^\prime(u_\alpha)\bigr) \right),
$$
or, more concisely,
$$
x^+-u_\alpha=(E-\mu\alpha\mathcal H)(x-u_\alpha)-\mu\alpha\bigl(J_\delta^\prime(u_\alpha)-J^\prime(u_\alpha)\bigr), \eqno{(31)}
$$
where
$$
\mathcal H=\int\limits_0^1 T_\alpha^{\delta\,\prime\prime}(u_\alpha+t(x-u_\alpha))dt.
$$

We next study the spectrum of the self-adjoint operator $\mathcal H$. To this end, we estimate the operators $T_\alpha^{\delta\,\prime\prime}(u_\alpha+t(x-u_\alpha))$ appearing under the integral sign. We use the standard ordering of self-adjoint operators: the notation $A\geq 0$ means that $A$ is positive semidefinite, while $A\geq B$ is equivalent to $A-B\geq 0$. In particular, $J^{\prime\prime}(x^*)\geq 0$. If $mE\leq A\leq ME$ for some numbers $m,M\in\mathbb{R}$, then the spectrum $\sigma(A)$ of the self-adjoint operator $A$ lies in the interval $[m,M]$. Recall that $E$ is the identity operator on $H$.

For any point $z$ on the line segment joining $u_\alpha$ and $x$, Theorem 3.1 and assumption (30) give
$$
\|z-x^*\|\leq\|z-u_\alpha\|+\|u_\alpha-x^*\|\leq(\varkappa+\widehat C_0)\alpha, \eqno{(32)}
$$
and therefore
\begin{equation*}
\begin{split}
&T_\alpha^{\prime\prime}(z)=J^{\prime\prime}(z)+2\alpha E=J^{\prime\prime}(x^*)+(J^{\prime\prime}(z)-J^{\prime\prime}(x^*))+2\alpha E\geq\\
&\geq \bigl(2\alpha-\|J^{\prime\prime}(z)-J^{\prime\prime}(x^*)\|_{L(H)}\bigr)E\geq\\
&\geq (2\alpha-L\|z-x^*\|)E\geq(2-L(\varkappa+\widehat C_0))\alpha E;
\end{split}
\tag{33}
\end{equation*}
\begin{equation*}
\begin{split}
T_\alpha^{\prime\prime}(z)&\leq\|T_\alpha^{\prime\prime}(z)\|_{L(H)}E\leq (\|J^{\prime\prime}(x^*)\|_{L(H)}+\|J^{\prime\prime}(z)-J^{\prime\prime}(x^*)\|_{L(H)}+2\alpha) E\leq\\
&\leq \bigl(\|J^{\prime\prime}(x^*)\|_{L(H)}+L(\varkappa+\widehat C_0)\alpha+2\alpha\bigr) E.
\end{split}
\end{equation*}
Define
$$
C_4=2-L(\varkappa+\widehat C_0), \eqno{(34)}
$$
and assume that $C_4>0$, or equivalently,
$$
\varkappa<\frac{2}{L}-\widehat C_0. \eqno{(35)}
$$
Such a choice of $\varkappa$ is possible by (24). We have
$$
C_4\alpha E\leq T_\alpha^{\prime\prime}(z)\leq \bigl(\|J^{\prime\prime}(x^*)\|_{L(H)}+L(\varkappa+\widehat C_0)\alpha+2\alpha\bigr) E.
$$
Combining this two-sided estimate with (27) and (32), we obtain
\begin{equation*}
\begin{split}
&T_\alpha^{\delta\,\prime\prime}(z)=J_\delta^{\prime\prime}(z)+2\alpha E=J_\delta^{\prime\prime}(z)-J^{\prime\prime}(z)+T_\alpha^{\prime\prime}(z)\geq\\
&\geq (-\delta(1+\|z-\xi\|)+C_4\alpha)E\geq \bigl(C_4\alpha-\delta(1+\|x^*-\xi\|+(\varkappa+\widehat C_0)\alpha)\bigr)E;
\end{split}
\end{equation*}
\begin{equation*}
\begin{split}
&T_\alpha^{\delta\,\prime\prime}(z)\leq\|T_\alpha^{\delta\,\prime\prime}(z)\|_{L(H)}E\leq \bigl(\|J_\delta^{\prime\prime}(z)-J^{\prime\prime}(z)\|_{L(H)}+\|T_\alpha^{\prime\prime}(z)\|_{L(H)}\bigr)E\leq\\
&\leq \bigl(\delta(1+\|x^*-\xi\|+(\varkappa+\widehat C_0)\alpha)+ \|J^{\prime\prime}(x^*)\|_{L(H)}+L(\varkappa+\widehat C_0)\alpha+2\alpha \bigr)E.
\end{split}
\end{equation*}
Set
$$
C_5=\|J^{\prime\prime}(x^*)\|_{L(H)}+\frac{C_4\bar\alpha}{2}+L(\varkappa+\widehat C_0)\bar\alpha+2\bar\alpha \eqno{(36)}
$$
and require that
$$
0\leq \delta\leq \frac{C_4\alpha}{2(1+\|x^*-\xi\|+(\varkappa+\widehat C_0)\bar\alpha)}. \eqno{(37)}
$$
Then
$$
\frac{C_4\alpha}{2}E\leq T_\alpha^{\delta\,\prime\prime}(z)\leq C_5E.
$$
Integrating these operator inequalities along the line segment joining $u_\alpha$ and $x$, we obtain
$$
\frac{C_4\alpha}{2}E\leq \mathcal H\leq C_5E,
$$
i.e., $\sigma(\mathcal H)\subset[C_4\alpha/2,C_5]$.

We now return to (31) and estimate the norm of the operator $E-\mu\alpha\mathcal H$ on its right-hand side:
$$
\|E-\mu\alpha\mathcal H\|_{L(H)}\leq \max\limits_{\lambda\in[C_4\alpha/2,C_5]}|1-\mu\alpha\lambda|.
$$
Clearly, $1-\mu\alpha\lambda\leq 1-C_4\mu\alpha^2/2$ for $\lambda\in[C_4\alpha/2,C_5]$. Suppose that
$$
\mu\alpha\left(C_5+\frac{C_4\alpha}{2}\right)\leq 2. \eqno{(38)}
$$
Then we also have $-(1-\mu\alpha\lambda)\leq C_5\mu\alpha-1\leq 1-C_4\mu\alpha^2/2$. Consequently,
$$
\|E-\mu\alpha\mathcal H\|_{L(H)}\leq 1-\frac{C_4\mu\alpha^2}{2}.
$$
Note that the nonnegativity of the right-hand side of this inequality follows from (36) and (38). Furthermore, condition (8) and Theorem 3.1 give
$$
\|J_\delta^\prime(u_\alpha)-J^\prime(u_\alpha)\|\leq \delta(1+\|x^*-\xi\|+\widehat C_0\alpha).
$$
Now representation (31) yields
$$
\|x^+-u_\alpha\|\leq\|E-\mu\alpha\mathcal H\|_{L(H)}\|x-u_\alpha\|+\mu\alpha\|J_\delta^\prime(u_\alpha)-J^\prime(u_\alpha)\|\leq
$$
$$
\leq \left(1-\frac{C_4\mu\alpha^2}{2}\right)\|x-u_\alpha\|+\mu\alpha\delta(1+\|x^*-\xi\|+\widehat C_0\alpha).
$$
Set
$$
C_6=1+\|x^*-\xi\|+\widehat C_0\bar\alpha.
$$
Then
$$
\|x^+-u_\alpha\|\leq\left(1-\frac{C_4\mu\alpha^2}{2}\right)\|x-u_\alpha\|+C_6\mu\alpha\delta. \eqno{(39)}
$$
We have proved the following lemma.

\medskip

{\bf Lemma 4.1.} Suppose that the assumptions of Theorem 3.1 hold, that the parameter $\alpha$ satisfies (29), and that (30) holds with $\varkappa>0$ satisfying (35). Suppose, in addition, that inequality (38) holds with the constants $C_4$ and $C_5$ defined in (34) and (36). Finally, suppose that the error level $\delta$ satisfies (37). Then estimate (39) holds for one step (28) of the iteratively regularized gradient method with an approximate functional $J_\delta$ satisfying conditions (8) and (27).

\medskip

{\bf Remark 4.1.} Inequality (33), together with the condition $C_4>0$, implies that the functional $T_\alpha$ is strongly convex on the ball defined by (30). This, in turn, implies the Fej\'er-type property
$$
(T_\alpha^\prime(x),x-u_\alpha)\geq C_4\alpha\|x-u_\alpha\|^2, \eqno{(40)}
$$
which could also be used to prove an estimate of the form (39). In particular, if $\delta=0$, then
$$
\|x^+-u_\alpha\|^2=\|x-u_\alpha-\mu\alpha T_\alpha^\prime(x)\|^2=
$$
$$
=\|x-u_\alpha\|^2-2\mu\alpha(T_\alpha^\prime(x),x-u_\alpha)+\mu^2\alpha^2\|T_\alpha^\prime(x)\|^2=
$$
$$
=\|x-u_\alpha\|^2-2\mu\alpha(T_\alpha^\prime(x),x-u_\alpha)+\mu^2\alpha^2\|T_\alpha^\prime(x)-T_\alpha^\prime(u_\alpha)\|^2.
$$
By (40), the second term on the right-hand side, as well as the other two terms, can be estimated in terms of $\|x-u_\alpha\|^2$. Such an argument is described in [6] and [11, Chap.~2, \S~8] for iterative methods for solving operator equations. In Lemma 4.1, we used an alternative approach based on studying the spectrum of $\mathcal H$ instead of applying (40). The Fej\'er-based argument also yields an error estimate of the same order. However, it is more cumbersome, especially when $\delta>0$, and is less convenient for deriving explicit formulas for the constants.

\medskip

For the subsequent arguments, we need to strengthen condition (35) on the choice of $\varkappa$. Namely, we require that
$$
\widehat C_0<\varkappa<\frac{2}{L}-\widehat C_0. \eqno{(41)}
$$
A constructive choice of $\varkappa$ is discussed at the end of this section. Condition (41) can be satisfied if and only if
$$
L\widehat C_0<1, \eqno{(42)}
$$
which is stronger than condition (24) used above.

Let us determine when (42) holds. If we set $t=L\|w\|$, then $t\in[0,1)$ by condition (12), and
$$
L\widehat C_0=\frac{t}{1-t}\left(t+\sqrt{t^2+(1-t)\left(4+\frac{8}{3}t\right)}\right). \eqno{(43)}
$$
We arrive at the inequality
$$
\frac{t}{1-t}\left(t+\sqrt{t^2+(1-t)\left(4+\frac{8}{3}t\right)}\right)<1,\quad t\in[0,1). \eqno{(44)}
$$
Equivalently, it can be written as
$$
t\sqrt{t^2+(1-t)\left(4+\frac{8}{3}t\right)}<1-t-t^2,\quad t\in(0,1).
$$
The left-hand side is always nonnegative, whereas the right-hand side is positive only for $t\in\left[0,\frac{\sqrt{5}-1}{2}\right)$. Squaring both sides yields the fourth-degree inequality
$$
\frac{8}{3}t^4+\frac{10}{3}t^3-5t^2-2t+1>0,\quad t\in\left[0,\frac{\sqrt{5}-1}{2}\right).
$$
It can be reduced to a cubic inequality by factoring the left-hand side:
$$
(t-1)\left(\frac{8}{3}t^3+6t^2+t-1\right)>0;
$$
$$
8t^3+18t^2+3t-3<0,\quad t\in\left[0,\frac{\sqrt{5}-1}{2}\right).
$$
The cubic polynomial on the left-hand side is increasing for $t\geq 0$ and has a unique positive root $t_0\approx 0.316$ belonging to the interval $\left[0,\frac{\sqrt{5}-1}{2}\right)$. Thus, inequality (44) holds for $t\in[0,t_0)$, and hence, in view of (12), condition (42) is equivalent to
$$
L\|w\|<t_0,\,\,{\textrm{where}}\,\,8t_0^3+18t_0^2+3t_0-3=0,\,t_0>0\,\,(t_0\approx 0.316). \eqno{(45)}
$$
Note that condition (12) follows automatically from (45). We have proved the following lemma.

\medskip

{\bf Lemma 4.2.} Suppose that conditions (2), (10), and (45) hold. Then all the assumptions of Theorem 3.1 are satisfied, and there exist values of $\varkappa$ satisfying (41).

\medskip

We now proceed to prove the main result of this section.

\medskip

{\bf Theorem 4.1.} Suppose that conditions (2), (10), and (45) hold, and that the approximate functional $J_\delta$ is twice continuously Fr\'echet differentiable and satisfies inequalities (8) and (27). If the functional $J$ is known exactly, set $J_\delta=J$, $\delta=0$. Further, let the constant $\varkappa$ satisfy (41), which is always possible by Lemma 4.2, and fix $\bar\alpha>0$, $\mu>0$, $\alpha_0>0$, $0<s<2/5$. Consider the iterative process
$$
x_{n+1}=x_n-\mu\alpha_n\bigl(J_\delta^\prime(x_n)+2\alpha_n(x_n-\xi)\bigr),\quad
\alpha_n=\alpha_0(n+k)^{-s},\quad n\in\mathbb{N}. \eqno{(46)}
$$
Choose $k$ in (46) sufficiently large, independently of $\delta$, and suppose that
$$
\|x_1-x^*\|\leq(\varkappa-\widehat C_0)\alpha_1, \eqno{(47)}
$$
where the right-hand side is positive by (41). If $\delta=0$, then the error estimate
$$
\|x_n-x^*\|\leq(\varkappa+\widehat C_0)\alpha_n,\quad n\in\mathbb{N} \eqno{(48)}
$$
holds. If $\delta>0$, suppose that the stopping index is defined by
$$
N=\min\{n\geq 1\,|\,\alpha_n\leq K_1\sqrt\delta\}, \eqno{(49)}
$$
where
$$
K_1\geq\sqrt{\max\left\{\frac{2C_7\bar\alpha}{C_4},\frac{4C_6}{\varkappa C_4}\right\}},\quad C_7=1+\|x^*-\xi\|+(\varkappa+\widehat C_0)\bar\alpha. \eqno{(50)}
$$
Then
$$
\|x_N-x^*\|\leq(\varkappa+\widehat C_0) K_1\sqrt\delta. \eqno{(51)}
$$

{\bf Proof.} We first ensure that the conditions of Lemma 4.1 are satisfied. Since $\alpha_1=\alpha_0(k+1)^{-s}\to 0$ as $k\to\infty$, we can choose $k$ sufficiently large so that inequalities (29) and (38) hold simultaneously for $\alpha=\alpha_1$. We make this choice throughout the remainder of the proof. The sequence $\{\alpha_n\}$ is decreasing, and therefore the same two inequalities hold for all $\alpha=\alpha_n$, $n\in\mathbb{N}$. In particular, all $\alpha_n\leq\bar\alpha$. Further, if $\delta=0$, condition (37) holds for every $n\in\mathbb{N}$. If $\delta>0$, then for $n<N$,
$$
\delta<\frac{\alpha_n^2}{K_1^2}\leq\frac{C_4\alpha_n^2}{2C_7\bar\alpha}\leq\frac{C_4\alpha_n}{2C_7}.
$$
This is exactly condition (37) with $\alpha=\alpha_n$. Thus, all the parameter conditions of Lemma 4.1 hold for $n<N$.

Lemma 4.1 applies when the point $x$ considered there belongs to the ball (30) centered at $u_\alpha$. For $\alpha=\alpha_n$, $n\in\mathbb{N}$, the corresponding balls form the tube $\bigcup\limits_{n=1}^\infty\overline B(u_{\alpha_n},\varkappa\alpha_n)$. Here, $\overline B(x,r)$ is the closed ball centered at $x\in H$ with radius $r>0$. Set
$$
d_n=\|x_n-u_{\alpha_n}\|.
$$
Then the condition for the $n$th ball to contain the point $x_n$ can be written as
$$
d_n\leq\varkappa\alpha_n. \eqno{(52)}
$$
We prove by induction that $x_n$ belongs to the $n$th ball for every $n\in\mathbb{N}$ when $\delta=0$, and for every $n\leq N$ when $\delta>0$. Thus, until the iterations are stopped, all the approximations obtained belong to the indicated tube.

The base case follows from (47) and Theorem 3.1:
$$
d_1\leq\|x_1-x^*\|+\|u_{\alpha_1}-x^*\|\leq\varkappa\alpha_1.
$$
Assume that $d_n\leq\varkappa\alpha_n$ and, if $\delta>0$, that $n<N$. Applying Lemmas 4.1 and 3.1 gives
\begin{equation*}
\begin{split}
&d_{n+1}\leq\|x_{n+1}-u_{\alpha_n}\|+\|u_{\alpha_n}-u_{\alpha_{n+1}}\|\leq\\
&\leq\left(1-\frac{C_4\mu\alpha_n^2}{2}\right)d_n+C_6\mu\alpha_n\delta+C_3\frac{\Delta_n}{\sqrt{\alpha_n}}.
\end{split}
\tag{53}
\end{equation*}
Here we have set
$$
\Delta_n=\alpha_n-\alpha_{n+1}.
$$
Substituting the induction hypothesis $d_n\leq\varkappa\alpha_n$ into (53) shows that $d_{n+1}\leq\varkappa\alpha_{n+1}$ follows if
$$
C_6\mu\alpha_n\delta+\left(\frac{C_3}{\sqrt{\alpha_n}}+\varkappa\right)\Delta_n\leq\frac{C_4\mu\varkappa}{2}\alpha_n^3.
$$
By (49) and (50), if $\delta>0$ and $n<N$, then
$$
C_6\mu\alpha_n\delta<\frac{C_6\mu}{K_1^2}\alpha_n^3\leq\frac{C_4\mu\varkappa}{4}\alpha_n^3,
$$
and this inequality also holds when $\delta=0$. Therefore, the induction step is valid if
$$
\Delta_n\leq C_8\alpha_n^{7/2},\,\,{\textrm{where}}\,\,C_8=\frac{\varkappa C_4\mu}{4(C_3+\varkappa\sqrt{\bar\alpha})}>0. \eqno{(54)}
$$
Let us verify (54). We have
$$
\Delta_n=\alpha_0\bigl((n+k)^{-s}-(n+k+1)^{-s}\bigr)\leq s\alpha_0(n+k)^{-s-1},
$$
and hence
$$
\frac{\Delta_n}{\alpha_n^{7/2}}\leq \frac{s\alpha_0(n+k)^{-s-1}}{\alpha_0^{7/2}(n+k)^{-7s/2}}=s\alpha_0^{-5/2}(n+k)^{5s/2-1}.
$$
Since $s<2/5$ by the hypothesis of the theorem, the right-hand side decreases as $n$ increases. Hence, condition (54) holds for every $n\geq 1$ if $k$ is additionally chosen so that
$$
k+1\geq\left(\frac{s}{C_8\alpha_0^{5/2}}\right)^{\frac{1}{1-5s/2}}. \eqno{(55)}
$$
We can increase the value of $k$ chosen at the beginning of the proof until (55) is satisfied, without violating any of the conditions. We emphasize that condition (47) in Theorem 4.1 is imposed only after the final choice of $k$.

Thus, the induction proves (52) for all $n\in\mathbb{N}$ when $\delta=0$ and for $n\leq N$ when $\delta>0$. By Lemma 4.2, all the assumptions of Theorem 3.1 hold, and therefore
$$
\|x_n-x^*\|\leq d_n+\|u_{\alpha_n}-x^*\|\leq(\varkappa+\widehat C_0)\alpha_n.
$$
This proves estimate (48) for $\delta=0$. For $\delta>0$, setting $n=N$ and using (49) also gives the required estimate (51). {\bf Theorem proved.}

\medskip

{\bf Remark 4.2.} It is readily seen that the endpoint value $s=2/5$ is also admissible in Theorem 4.1 if, instead of increasing $k$ according to (55), one imposes the additional condition
$$
\frac{2}{5}\alpha_0^{-5/2}\leq C_8.
$$

\medskip\medskip

We conclude with some remarks on the constructive choice of the parameters in Theorem 4.1. When method (46) is applied to specific optimization problems, the constant $L$ can be expected to be known, whereas $\|w\|$ need not be known. In any event, most results concerning error estimates for methods of solving ill-posed problems under a source condition on the sought solution assume that $\|w\|$ is sufficiently small. We therefore assume that a bound $\|w\|\leq w_0$ is known for some $w_0>0$ such that condition (45) can be verified. This condition, in turn, implies inequality (42), which shows that the choice $\varkappa=1/L$ certainly satisfies (41), even if $\widehat C_0$ is unknown. On the other hand, knowledge of $w_0$ makes it possible to use formula (43) to obtain an upper bound for $\widehat C_0$, thereby enlarging the admissible range of $\varkappa$ and refining condition (47) on the proximity of $x_1$ to $x^*$. Note that the right-hand side of (43) tends to zero as $t\to 0$; therefore, as $w_0\to 0$, the guaranteed upper bound for $\widehat C_0$ also tends to zero. Further, choosing the constant $K_1$ in stopping rule (49) requires upper bounds for $C_6$ and $C_7$ and a lower bound for $C_4$. All three bounds can be obtained if upper bounds for $\|x^*-\xi\|$ and $\widehat C_0$ are known.

The rule for choosing the parameter $k$ in the iterative process (46) is also specified in the proof of the theorem: conditions (29) and (38) must hold for $\alpha=\alpha_1=\alpha_0(k+1)^{-s}$, as must condition (55), with $\alpha_0$, $s$, $\varkappa$, and $\mu$ already fixed. This requires upper and lower bounds for the constant $C_4$, as well as upper bounds for the constants $C_3$ and $C_5$. Analysis of the explicit formulas above for these constants shows that an upper bound for $\|J^{\prime\prime}(x^*)\|_{L(H)}$ is also required. All these requirements are feasible in principle.

We also note that $\bar\alpha$ in the assumptions of Theorem 4.1 need not be fixed in advance and can instead be set equal to $\alpha_1=\alpha_0(k+1)^{-s}$. Condition (29) then holds automatically for all $\alpha_n$. In this case, the constants $C_3$, $C_5$, $C_6$, $C_7$, and $C_8$ depend on $\alpha_1$ and hence on $k$, while the choice of $k$ depends on the values of these constants. However, this does not create a logical circle. Indeed, as $k$ increases, $\alpha_1$ decreases; consequently, $C_3$ and $C_5$ also decrease, while $C_8$ increases. Thus, $k$ can still be chosen sufficiently large for conditions (38) and (55) to hold. Moreover, $C_6$ and $C_7$ decrease as $\bar\alpha$ decreases, which weakens the condition on $K_1$ in (50) and may improve estimate (51). Therefore, the choice $\bar\alpha=\alpha_1$ is optimal, since, by (29), this is the smallest possible value of $\bar\alpha$.

\section{An error estimate for the iteratively regularized gradient method: second proof}

We now present the second proof of the error estimate for the iteratively regularized gradient method. This proof uses neither structural properties of the Tikhonov functional nor its minimizers. Moreover, we impose only condition (8) on the approximate functional $J_\delta$, without the additional requirement (27). The smoothness condition (2) and source condition (10) are still assumed to hold, while instead of (45) we use the slightly stronger condition
$$
L\|w\|<\frac{1}{4}. \eqno{(56)}
$$
Fix $\mu>0$, $\alpha_0>0$, $0<s<1/2$, and consider the iterative process (46). Conditions on $k$ will be derived below.

By (2) and the identity $J^\prime(x^*)=0_H$, we have
$$
J^\prime(x)=J^{\prime\prime}(x^*)(x-x^*)+R(x),\quad \|R(x)\|\leq\frac{L}{2}\|x-x^*\|^2,\quad x\in H. \eqno{(57)}
$$
Combining this representation with (8) and (10), we obtain
\begin{equation*}
\begin{split}
&T_{\alpha_n}^{\delta\,\prime}(x_n)=J_\delta^\prime(x_n)+2\alpha_n(x_n-\xi)=(J_\delta^\prime(x_n)-J^\prime(x_n))+\\
&+(J^{\prime\prime}(x^*)(x_n-x^*)+R(x_n))+2\alpha_n(x_n-x^*)+2\alpha_nJ^{\prime\prime}(x^*)w=\\
&=(J_\delta^\prime(x_n)-J^\prime(x_n))+J^{\prime\prime}(x^*)z_n+2\alpha_nz_n-4\alpha_n^2w+R(x_n),
\end{split}
\end{equation*}
where we have set
$$
z_n=x_n-x^*+2\alpha_nw. \eqno{(58)}
$$
Combining this representation with (46), we obtain
\begin{equation*}
\begin{split}
&z_{n+1}=z_n-2(\alpha_n-\alpha_{n+1})w-\\
&-\mu\alpha_n\bigl((J_\delta^\prime(x_n)-J^\prime(x_n))+J^{\prime\prime}(x^*)z_n+2\alpha_nz_n-4\alpha_n^2 w+R(x_n)\bigr)=\\
&=Q_nz_n-\mu\alpha_n\bigl(J_\delta^\prime(x_n)-J^\prime(x_n)\bigr)-\mu\alpha_nR(x_n)+4\mu\alpha_n^3 w-2(\alpha_n-\alpha_{n+1})w,
\end{split}
\tag{59}
\end{equation*}
where
$$
Q_n=(1-2\mu\alpha_n^2)E-\mu\alpha_n J^{\prime\prime}(x^*).
$$

Suppose that
$$
\|z_1\|\leq q\alpha_1 \eqno{(60)}
$$
with some constant $q>0$. Inequality (60) will serve as the base case of an induction. We next determine conditions under which the corresponding inequality
$$
\|z_n\|\leq q\alpha_n \eqno{(61)}
$$
can be proved by induction for all $n\in\mathbb{N}$, or at least for $n\leq N$ with a suitable stopping index $N$.

Assume inductively that (61) holds for some $n\in\mathbb{N}$. We show that $\|z_{n+1}\|\leq q\alpha_{n+1}$. By (8) and (10),
$$
\|J_\delta^\prime(x_n)-J^\prime(x_n)\|\leq \delta(1+\|x_n-x^*\|+\|J^{\prime\prime}(x^*)w\|)\leq
$$
$$
\leq \delta(1+\|z_n-2\alpha_n w\|+\|J^{\prime\prime}(x^*)w\|)\leq \delta(1+\|J^{\prime\prime}(x^*)w\|+(q+2\|w\|)\alpha_n).
$$
Set
$$
C_9=1+\|J^{\prime\prime}(x^*)w\|+(q+2\|w\|)\alpha_1.
$$
Using $\alpha_n\leq\alpha_1$, we arrive at
$$
\|J_\delta^\prime(x_n)-J^\prime(x_n)\|\leq C_9\delta. \eqno{(62)}
$$

We now estimate the norm of the operator $Q_n$ in (59). The spectrum of $J^{\prime\prime}(x^*)$ satisfies
$$
\sigma(J^{\prime\prime}(x^*))\subset\left[0,\|J^{\prime\prime}(x^*)\|_{L(H)}\right],
$$
and therefore
$$
\sigma(Q_n)\subset\left[1-2\mu\alpha_n^2-\mu\alpha_n\|J^{\prime\prime}(x^*)\|_{L(H)},\,1-2\mu\alpha_n^2\right].
$$
It follows that the estimate
$$
\|Q_n\|_{L(H)}\leq 1-2\mu\alpha_n^2 \eqno{(63)}
$$
holds if we require
$$
1-2\mu\alpha_n^2-\mu\alpha_n\|J^{\prime\prime}(x^*)\|_{L(H)}\geq -(1-2\mu\alpha_n^2),
$$
equivalently, $\|J^{\prime\prime}(x^*)\|_{L(H)}\leq 2(1-2\mu\alpha_n^2)/(\mu\alpha_n)=g(\alpha_n)$, where
$$
g(\alpha)=\frac{2(1-2\mu\alpha^2)}{\mu\alpha}.
$$
The function $g$ is decreasing for $\alpha>0$; therefore, since $\alpha_n\leq\alpha_1$, it suffices to require $\|J^{\prime\prime}(x^*)\|_{L(H)}\leq g(\alpha_1)$, i.e.,
$$
\|J^{\prime\prime}(x^*)\|_{L(H)}\leq \frac{2(1-2\mu\alpha_1^2)}{\mu\alpha_1}. \eqno{(64)}
$$
Note that $g(\alpha)\to +\infty$ as $\alpha\to 0$, and that as $k\to\infty$ we have $\alpha_1=\alpha_0(k+1)^{-s}\to 0$. Thus, by choosing $k$ sufficiently large, we can make $g(\alpha_1)$ arbitrarily large and, in particular, ensure that it exceeds $\|J^{\prime\prime}(x^*)\|_{L(H)}$. We therefore choose $k\in\mathbb{N}$ sufficiently large so that (64) holds. Then estimate (63) holds for every $n\in\mathbb{N}$.

We now apply estimates (57), (62), and (63), the induction hypothesis (61), and the following consequence of (58) to (59):
$$
\|x_n-x^*\|\leq \|z_n\|+2\alpha_n\|w\|\leq (q+2\|w\|)\alpha_n.
$$
We obtain
$$
\|z_{n+1}\|\leq (1-2\mu\alpha_n^2)q\alpha_n+C_9\mu\alpha_n\delta+
$$
$$
+\frac{\mu L}{2}(q+2\|w\|)^2\alpha_n^3+4\mu\|w\|\alpha_n^3+2\|w\|(\alpha_n-\alpha_{n+1}).
$$
We rewrite this estimate as
$$
\|z_{n+1}\|\leq q\alpha_n-\mu\Phi(q)\alpha_n^3+C_9\mu\alpha_n\delta+2\|w\|(\alpha_n-\alpha_{n+1}), \eqno{(65)}
$$
where
$$
\Phi(q)=2q-4\|w\|-\frac{L}{2}(q+2\|w\|)^2. \eqno{(66)}
$$
An elementary analysis shows that $\Phi(q)$ attains its maximum value $\frac{2}{L}-8\|w\|$ at $q=\frac{2}{L}-2\|w\|$. Consequently, condition (56) is necessary and sufficient for the existence of a $q>0$ such that $\Phi(q)>0$. We henceforth assume that
$$
2q-4\|w\|-\frac{L}{2}(q+2\|w\|)^2>0 \eqno{(67)}
$$
holds; equivalently,
$$
\frac{2}{L}\left(1-L\|w\|-\sqrt{1-4L\|w\|}\right)<q<\frac{2}{L}\left(1-L\|w\|+\sqrt{1-4L\|w\|}\right).
$$

A sufficient condition for the right-hand side of (65) not to exceed $q\alpha_{n+1}$ is
$$
C_9\mu\alpha_n\delta+(q+2\|w\|)(\alpha_n-\alpha_{n+1})\leq\mu\Phi(q)\alpha_n^3. \eqno{(68)}
$$
Observe that
$$
\frac{\alpha_n-\alpha_{n+1}}{\alpha_n^3}\leq \frac{s\alpha_0(n+k)^{-s-1}}{\alpha_0^3 (n+k)^{-3s}}=s\alpha_0^{-2}(n+k)^{2s-1}. \eqno{(69)}
$$
If $s<1/2$, the right-hand side decreases with $n$. Fix $\theta\in(0,1)$ and increase $k$ sufficiently so that
$$
k+1\geq\left(\frac{s(q+2\|w\|)}{(1-\theta)\mu\Phi(q)\alpha_0^2}\right)^{\frac{1}{1-2s}}. \eqno{(70)}
$$
Note that increasing $k$ in this way does not invalidate condition (64), which has already been ensured. Then $n+k$ is also greater than or equal to the right-hand side of (70), and (69) implies
$$
\frac{\alpha_n-\alpha_{n+1}}{\alpha_n^3}\leq\frac{(1-\theta)\mu\Phi(q)}{q+2\|w\|}.
$$
Hence, for (68) to hold, it suffices that
$$
\delta\leq\frac{\theta\Phi(q)}{C_9}\alpha_n^2. \eqno{(71)}
$$

If $\delta=0$, inequality (71) holds, and hence estimate (61) is valid for every $n\in\mathbb{N}$ provided that the base case (60) holds. Therefore, the estimate
$$
\|x_1-x^*\|\leq (q-2\|w\|)\alpha_1 \eqno{(72)}
$$
implies
$$
\|x_n-x^*\|\leq(q+2\|w\|)\alpha_n,\quad n\in\mathbb{N}. \eqno{(73)}
$$
Note that $q>2\|w\|$ here by (67).

Now let $\delta>0$. Choose an arbitrary constant
$$
K_2\geq\sqrt{\frac{C_9}{\theta\Phi(q)}}, \eqno{(74)}
$$
and let
$$
N=\min\{n\geq1\,|\,\alpha_n\leq K_2\sqrt\delta\}. \eqno{(75)}
$$
Then condition (71) holds for every $n<N$, and therefore the base case (60) implies estimate (61) for every $n\leq N$. Hence, if (72) holds, then
$$
\|x_n-x^*\|\leq(q+2\|w\|)\alpha_n,\quad 1\leq n\leq N.
$$
Setting $n=N$ here and using (75), we arrive at the final estimate
$$
\|x_N-x^*\|\leq(q+2\|w\|) K_2\sqrt\delta. \eqno{(76)}
$$
Thus, we have proved the following theorem.

\medskip

{\bf Theorem 5.1.} Suppose that conditions (2), (8), (10), and (56) hold. Fix $0<\theta<1$, $\mu>0$, $\alpha_0>0$, $0<s<1/2$, and choose $q>0$ so that (67) holds. Assume that $k$ in (46) is chosen sufficiently large, independently of $\delta$, and that condition (72) holds; its right-hand side is positive by (67). Then estimate (73) holds when $\delta=0$, while estimate (76) holds when $\delta>0$ and the stopping index is chosen according to rule (75). The constant $K_2$ in (75) and (76) must satisfy condition (74), and the expression $\Phi(q)$ is defined in (66).

\medskip

{\bf Remark 5.1.} The endpoint value $s=1/2$ is also admissible in Theorem 5.1, provided that condition (70) is replaced by
$$
\frac{1}{2\alpha_0^2}\leq\frac{\mu(1-\theta)\Phi(q)}{q+2\|w\|}.
$$

We conclude with a brief discussion of the constructive choice of the parameters in Theorem 5.1. Knowledge of a bound $\|w\|\leq w_0$ with $Lw_0<1/4$ guarantees condition (56) and makes it possible to choose $q$ so that (67) holds. More precisely, since the left-hand side of (67) decreases as $\|w\|$ increases, it suffices to require $\Phi_0=2q-4w_0-\frac{L}{2}(q+2w_0)^2>0$, which can certainly be satisfied because $Lw_0<1/4$. Then $\Phi(q)\geq \Phi_0$. Together with an upper bound for $\|J^{\prime\prime}(x^*)\|_{L(H)}$, these estimates allow one to choose $k$ so that (64) and (70) hold, obtain an upper bound for $C_9$, and then choose $K_2$ according to (74).

\section{Comparison of the two proofs}

Let us compare Theorems 4.1 and 5.1. They give error estimates of the same order for the iteratively regularized gradient method (46) applied to problem (1) of minimizing an arbitrary smooth functional $J$ on a Hilbert space. Specifically, the two theorems yield an $O(\alpha_n)$ estimate in the exact case and an $O(\sqrt\delta)$ estimate in the perturbed case with an a priori stopping rule. Theorem 5.1 has two significant advantages over Theorem 4.1: it does not require condition (27) on $J_\delta$, and it allows the larger range $0<s<1/2$ rather than $0<s<2/5$ in (46). On the other hand, condition (45) in Theorem 4.1 is weaker than the analogous condition (56) in Theorem 5.1.

Moreover, Theorem 4.1 may sometimes yield an error estimate with a smaller numerical coefficient than Theorem 5.1. As an example, take $\delta=0$, $L\|w\|=\frac{1}{4}-\varepsilon_1$, and $\varkappa=\widehat C_0+\varepsilon_2$ in Theorem 4.1, where $\varepsilon_1,\varepsilon_2>0$ are arbitrarily small. Then, by formula (43), as $\varepsilon_1\to 0$ we have $L\widehat C_0\to \frac{1+\sqrt{57}}{12}$. Denote this limit by $a$. Theorem 4.1 then gives the estimate $\|x_n-x^*\|\leq (2\widehat C_0+\varepsilon_2)\alpha_n$ with a numerical coefficient $2\widehat C_0+\varepsilon_2$ arbitrarily close to $2a/L\approx 1.425/L$. By contrast, for the same problem, Theorem 5.1 gives estimate (73), whose coefficient of $\alpha_n$ is $q+2\|w\|$. As $\varepsilon_1\to 0$, the admissible interval for $q$ shrinks to the maximizer $\frac{2}{L}-2\|w\|$ of $\Phi$. Therefore, as $\varepsilon_1\to 0$, the coefficient in the estimate of Theorem 5.1 is arbitrarily close to $2/L$ and thus exceeds the coefficient in the estimate of Theorem 4.1. If $\delta>0$, the coefficient $K_2$ in the estimate of Theorem 5.1 grows without bound as $\varepsilon_1\to 0$, making this estimate increasingly weaker, whereas the coefficient $K_1$ in the analogous estimate of Theorem 4.1 remains bounded.

Thus, neither theorem dominates the other in all respects.

\medskip

{\bf Remark 6.1.} Note that the limiting case $\varepsilon_1=\varepsilon_2=0$ in the example above is inadmissible because of conditions (41) and (56). As $\varepsilon_1\to 0$, Theorem 5.1 requires arbitrarily large values of $k$ in method (46) and of $K_2$ in the stopping rule, whereas, as $\varepsilon_2\to 0$, Theorem 4.1 requires an increasingly restrictive condition on the proximity of the initial point $x_1$ to $x^*$.

\section{The case of local smoothness}

In all the theorems above, the functional $J$ was assumed to be twice continuously Fr\'echet differentiable and to satisfy condition (2) on the entire space $H$. However, this condition was used mainly locally, along explicitly defined line segments. It can therefore be weakened with only minimal changes to the arguments. Namely, let $U\subset H$ be an open neighborhood of the minimizer $x^*$. Suppose that the functional $J$ is twice continuously Fr\'echet differentiable only on $U$, with
$$
\|J^{\prime\prime}(x)-J^{\prime\prime}(y)\|_{L(H)}\leq L\|x-y\|,\quad x,y\in U. \eqno{(77)}
$$
Outside $U$, the functional $J$ need not be defined at all. Conditions (8) and (27) on the approximate functional $J_\delta$ are also assumed to hold only on $U$. Source condition (10) is still assumed to hold, while the point $\xi$ need not belong to $U$.

We first derive a local version of Theorem 4.1. Suppose that condition (45) holds and that $\varkappa$ satisfies (41). Fix $\bar\alpha>0$ and require
$$
\overline B\bigl(x^*,(\varkappa+\widehat C_0)\bar\alpha\bigr)\subset U. \eqno{(78)}
$$
Recall that $\overline B(x,r)$ is the closed ball centered at $x\in H$ with radius $r>0$.

\medskip

{\bf Lemma 7.1.} Under conditions (10), (41), (45), (77), and (78), for every $0<\alpha\leq\bar\alpha$ there exists a unique point $u_\alpha$ that minimizes the functional $T_\alpha$ on the closed ball
$$
D_\alpha=\overline B(x^*,\varkappa\alpha).
$$
It lies in the interior of $D_\alpha$ and satisfies $T_\alpha^\prime(u_\alpha)=0_H$ and estimate (25). For any $0<\beta<\alpha\leq\bar\alpha$, estimate (26) also holds with the constant $C_3$ specified in Lemma 3.1.

{\bf Proof.} By (78), $D_\alpha\subset U$. Using (77), as in (33), we obtain
$$
T_\alpha^{\prime\prime}(x)\geq (2-L\varkappa)\alpha E,\quad x\in D_\alpha.
$$
Here $2-L\varkappa>0$ by (41). Thus, $T_\alpha$ is strongly convex and continuous on the closed convex set $D_\alpha$ and hence has a unique minimizer on this set [1, Chap.~1, \S~3]. We denote it by $u_\alpha$. This allows us to use the formulas from the preceding sections without changing notation. We emphasize, however, that this point need not be a global minimizer of $T_\alpha$.

Observe that formula (23) implies
$$
\widehat C_0=C_0(0)\geq \sqrt{C_2(0)}\geq 2\|w\|.
$$
By (41),
$$
\varkappa>\widehat C_0\geq 2\|w\|. \eqno{(79)}
$$
Hence, the point $x^*-2\alpha w$ belongs to $D_\alpha$, and therefore
$$
T_\alpha(u_\alpha)\leq T_\alpha(x^*-2\alpha w).
$$
This inequality is sufficient to derive estimate (14) in the proof of Theorem 2.1 with $\varepsilon=0$ and $z=x^*-2\alpha w$. Formulas (16), (21), and (22) remain valid because $[x^*,u_\alpha],\,[x^*,x^*-2\alpha w],\,[u_\alpha,u_\alpha+2\alpha w]\subset U$, where $[a,b]$ denotes the line segment with endpoints $a,b\in H$. In particular, the last inclusion follows from
$$
\|u_\alpha+2\alpha tw-x^*\|\leq\|u_\alpha-x^*\|+2\alpha\|w\|\leq (\varkappa+2\|w\|)\alpha\leq(\varkappa+\widehat C_0)\bar\alpha,\quad t\in[0,1].
$$
Here we have used estimate (79) and condition (78). Thus, both stages of the proof of Theorem 2.1 for $\varepsilon=0$ can be repeated without change and yield estimate (25).

Since $\widehat C_0<\varkappa$, the point $u_\alpha$ lies in the interior of $D_\alpha$. Consequently, it is a local minimizer of $T_\alpha$, and therefore $T_\alpha^\prime(u_\alpha)=0_H$. The proof of Lemma 3.1 uses only the identities $T_\alpha^\prime(u_\alpha)=T_\beta^\prime(u_\beta)=0_H$, estimate (25), the positive semidefiniteness of $J^{\prime\prime}(x^*)$, and the Lipschitz continuity of $J^{\prime\prime}$ along the line segment joining $u_\alpha$ and $u_\beta$, which is contained in $U$. Therefore, the proof of Lemma 3.1 and the value of the constant $C_3$ remain unchanged. {\bf Lemma proved.}

\medskip

We now prove a local version of Theorem 4.1.

\medskip

{\bf Theorem 7.1.} Suppose that the local condition (77) holds instead of the global condition (2), that the functional $J_\delta$ is twice continuously Fr\'echet differentiable on $U$, and that conditions (8) and (27) hold for $x\in U$. Suppose also that (10) and (45) hold, that $\varkappa$ satisfies (41), and that inclusion (78) holds for some $\bar\alpha>0$. Then, under the remaining assumptions and with the choice of parameters specified in Theorem 4.1, estimates (48) and (51) remain valid with the same constants and the same rules for choosing $k$ and the stopping index $N$.

{\bf Proof.} In the proof of Theorem 4.1, we replace the global minimizers of $T_\alpha$ by the points $u_\alpha$ constructed in Lemma 7.1. The proof in Section 4 uses only their stationarity and estimates (25) and (26), so all the necessary properties are preserved.

Indeed, if $\|x-u_\alpha\|\leq\varkappa\alpha$, $0<\alpha\leq\bar\alpha$, and $z$ lies on the line segment joining $u_\alpha$ and $x$, then inequality (32) holds, and therefore, by (78), this entire line segment lies in $U$. Hence, representation (31), operator estimates (33), and their consequences remain valid. The entire tube of balls constructed in the proof of Theorem 4.1 also lies in $U$. Conditions (8), (27), and (77) are applied only at points already shown to belong to $U$. This allows the proof of Theorem 4.1 to be extended without change to the case where $J$ is only locally smooth. {\bf Theorem proved.}

\medskip

We now prove a local version of Theorem 5.1.

\medskip

{\bf Theorem 7.2.} Suppose that the local condition (77) holds instead of (2), that the functional $J_\delta$ is continuously Fr\'echet differentiable on $U$, and that condition (8) holds for $x\in U$. Suppose also that (10) and (56) hold, that $q>0$ satisfies (67), and that the parameters $k$ and $x_1$ are chosen so that conditions (64), (70), and (72), as well as the inclusion
$$
\overline B\bigl(x^*,(q+2\|w\|)\alpha_1\bigr)
\subset U \eqno{(80)}
$$
hold. Then all the conclusions of Theorem 5.1, including estimates (73) and (76), remain valid with the same constants.

{\bf Proof.} Under condition (77), representation (57) holds for every point $x$ such that the line segment $[x^*,x]$ is contained in $U$. Condition (72) ensures the base case (60) and places the initial approximation $x_1$ in the ball appearing in (80). If the induction hypothesis (61) holds, then
$$
\|x_n-x^*\|\leq(q+2\|w\|)\alpha_n\leq(q+2\|w\|)\alpha_1.
$$
Consequently, the point $x_n$ and the entire line segment $[x^*,x_n]$ lie in the ball from (80). Therefore, the local versions of (57) and (8) apply to $x_n$. This validates the induction step and allows all the arguments in the proof of Theorem 5.1 to be used without change. {\bf Theorem proved.}

\medskip

Thus, the ball on which the smoothness of $J$ is required has radius $O(\bar\alpha)$ for the first proof and $O(\alpha_1)$ for the second. Neither the orders $O(\alpha_n)$ and $O(\sqrt\delta)$ of the resulting error estimates nor their coefficients change. We also note that, as in Section 4, $\bar\alpha$ need not be fixed in advance in Theorem 7.1; one can instead set $\bar\alpha=\alpha_1$.

\begin{center}
\textbf{References}
\end{center}

1. {\it Vasil'ev F.P.} Methods for Solving Extremal Problems. Moscow: Nauka, 1981 (in Russian).

2. {\it Kokurin M.Yu.} Theory of Regularization of Ill-Posed Problems. Moscow; Izhevsk: Institute of Computer Science, 2026 (in Russian).

3. {\it Kaltenbacher B., Neubauer A., Scherzer O.} Iterative Regularization Methods for Nonlinear Ill-Posed Problems. Berlin: Walter de Gruyter, 2008. DOI: 10.1515/9783110208276

4. {\it Kokurin M.M.} Improved Accuracy Estimation of the Tikhonov Method for Ill-Posed Optimization Problems in Hilbert Space. {\it Computational Mathematics and Mathematical Physics.} 2023. V.63. No.4. P.519--527. DOI: 10.1134/S0965542523040103

5. {\it Vasil'ev F.P.} Regularization of Ill-Posed Problems of Minimization in Approximately Specified Sets. {\it U.S.S.R. Computational Mathematics and Mathematical Physics.} 1980. V.20. No.1. P.41--53. DOI: 10.1016/0041-5553(80)90060-9

6. {\it Bakushinskii A.B.} An Iteratively Regularized Gradient Method for Solving Nonlinear Irregular Equations. {\it Computational Mathematics and Mathematical Physics.} 2004. V.44. No.5. P.759--765.

7. {\it Kokurin M.M.} Convergence Rate Estimates for an Iteratively Regularized Gradient Method with an a Posteriori Stopping Rule. {\it Russian Mathematics.} 2026. V.70. No.2. P.34--42. DOI: 10.3103/S1066369X26700131

8. {\it Kindermann S.} Convergence of the gradient method for ill-posed problems. {\it Inverse Problems and Imaging.} 2017. V.11. No.4. P.703--720. DOI: 10.3934/ipi.2017033

9. {\it Kaltenbacher B., Van Huynh K.} Iterative regularization for constrained minimization formulations of nonlinear inverse problems. {\it Computational Optimization and Applications.} 2022. V.81. P.569--611. DOI: 10.1007/s10589-021-00343-x

10. {\it Kokurin M.Yu.} Convergence Rate Estimates for Tikhonov's Scheme as Applied to Ill-Posed Nonconvex Optimization Problems. {\it Computational Mathematics and Mathematical Physics.} 2017. V.57. No.7. P.1101--1110. DOI: 10.1134/S0965542517070090

11. {\it Bakushinskii A.B., Kokurin M.Yu.} Algorithmic Analysis of Irregular Operator Equations. Moscow: LENAND, 2012 (in Russian).

12. {\it Vasin V.V.} Fundamentals of the Theory of Ill-Posed Problems. Novosibirsk: Siberian Branch of the Russian Academy of Sciences Publishing House, 2020 (in Russian).

\end{document}